\input graphicx
\input amssym.def
\input amssym
\magnification=1000
\baselineskip = 0.18truein
\lineskiplimit = 0.01truein
\lineskip = 0.01truein
\vsize = 8.7truein
\voffset = 0.1truein
\parskip = 0.10truein
\parindent = 0.3truein
\settabs 12 \columns
\hsize = 6.0truein
\hoffset = 0.1truein

\setbox\strutbox=\hbox{%
\vrule height .708\baselineskip
depth .292\baselineskip
width 0pt}
\font\caps=cmcsc10

\font\bigtenrm=cmr10 at 14pt

\def\sqr#1#2{{\vcenter{\vbox{\hrule height.#2pt
\hbox{\vrule width.#2pt height#1pt \kern#1pt
\vrule width.#2pt}
\hrule height.#2pt}}}}
\def\square{\mathchoice\sqr46\sqr46\sqr{3.1}6\sqr{2.3}4}

\centerline{\bigtenrm SOME CONDITIONALLY HARD PROBLEMS}
\centerline{\bigtenrm ON LINKS AND 3-MANIFOLDS}
\tenrm
\vskip 14pt
\centerline{MARC LACKENBY}
\vskip 18pt
\centerline{\caps 1. Introduction}
\vskip 6pt

Many decision problems in the theory of knots, links and 3-manifolds are known to be solvable. For example, the equivalence problem
for links in the 3-sphere was solved by Haken [12], Hemion [15] and Matveev [21]. Following the work of many mathematicians, including the proof of the
the Geometrisation Conjecture by Perelman [23, 24, 25], the homeomorphism problem for compact orientable 3-manifolds is now solved [4, 18, 26].
However, the complexity of these and other decision problems in low-dimensional topology remains poorly understood.
The problem of deciding whether a knot is the unknot is a good test case. Haken was the first to find an algorithmic solution to this problem [11].
It was shown to be in NP by Hass, Lagarias
and Pippenger [14], and in co-NP by work of Agol [1], Kuperberg [17] and the author [20]. However, no polynomial-time algorithm has yet been found.

Although increasingly good algorithms continue to be discovered for these and other topological problems, it seems likely that many will be
not solvable in polynomial time. It therefore is relevant to find lower bounds on the complexity of these problems. Unfortunately, because of
the current status of complexity theory, one must usually be content to show that a problem is hard, assuming some standard conjectures,
such as P $\not=$ NP. We term such problems as {\sl conditionally hard}. But it is perhaps surprising that there are few problems
in low-dimensional topology that are even known to be conditionally hard. Indeed, there are just four results along these lines, 
as far as the author is aware.
One is due to Agol, Hass and Thurston [2], who proved that the problem of deciding whether a knot in a compact orientable 3-manifold bounds a compact orientable
surface with genus at most $g$ is NP-complete. The second is due to Burton and Spreer [8], who showed that the problem of deciding whether an
ideal triangulation of a 3-manifold admits a `taut' structure is also NP-complete. The third is due to Burton, Colin de Verdi\`ere and de Mesmay [6], who
showed that the problem of deciding whether a vector in normal co-ordinate space for a triangulated 3-manifold is represented by
an immersed normal surface is NP-hard. The fourth is due to Burton, de Mesmay and Wagner [7], who showed that the problem
of determining whether a closed triangulated 
3-manifold contains an embedded non-orientable surface with a given genus is NP-hard. Notably, none of these problems is solely about 
knots and links in the 3-sphere. Indeed, the problem of deciding whether a knot in the 3-sphere bounds a compact orientable surface with 
genus $g$ is in NP $\cap$ co-NP [2, 20],
and hence by a standard conjecture in complexity theory (see page 95 of [10]), this problem is believed {\sl not} to be NP-complete.

It is the purpose of this article to establish that some natural problems on links and 3-manifolds are conditionally hard. Unlike previous work,
some of these are concerned solely with links in the 3-sphere. We call these {\sl classical links}. We now explain these problems in detail.

\noindent {\caps Upper bound for the Thurston complexity of an unoriented classical link}

The {\sl Thurston complexity} $\chi_-(S)$ of a compact orientable connected surface $S$ is defined to be $\max \{ 0, -\chi(S) \}$.
The Thurston complexity of a compact orientable surface $S$ with components $S_1, \dots, S_k$ is
defined to be $\sum_i \chi_-(S_i)$. 

\vfill\eject
A {\sl Seifert surface} for a classical link $L$ is a compact orientable surface $S$ embedded in $S^3$,
with no closed components, such that $\partial S = L$. The {\sl Thurston complexity} $\chi_-(L)$
is the minimal Thurston complexity for any Seifert surface for $L$. Crucially, we do not require that
an orientation on the link components is specified, and hence there is no requirement
that the orientation on the surface $S$ is compatible with a given orientation on the link $L$.
We also do not require that Seifert surfaces be connected, unlike some authors.
While connectedness is a useful requirement for Seifert surfaces in some contexts, it is
a hindrance here. Representatives of a second homology class in a 3-manifold that have minimal
Thurston complexity may necessarily be disconnected. (The exterior of the unlink provides an example
of this phenomenon). It therefore seems natural to permit Seifert surfaces to be disconnected.

In the decision problem {\caps Upper bound for the Thurston complexity of an unoriented classical link}, one is given a
diagram of a link $L$ in the 3-sphere and an integer $n$ in binary, and
one asks whether $\chi_-(L) \leq n$.

We will prove the following theorem, which is possibly our most definitive result.

\noindent {\bf Theorem 1.1.} {\sl The problem} {\caps Upper bound for the Thurston complexity of an unoriented classical link} {\sl is NP-complete.}

The method of proof owes much to Agol-Hass-Thurston [2] but new ideas are required to make their arguments work in the
setting of classical links.

\vskip 6pt
\noindent {\caps The homeomorphism problem for closed 3-manifolds}
\vskip 6pt

This fundamental problem takes, as its input, triangulations of two closed 3-manifolds, and asks whether these manifolds
are homeomorphic. It is shown to be conditionally hard, in the following sense.

\noindent {\bf Theorem 1.2.} {\caps The graph isomorphism problem} {\sl is Karp-reducible to}
{\caps the homeomorphism problem for closed 3-manifolds}.

For the definition of Karp-reducibility, see Definition 2.11 in [10]. Informally, this asserts that
{\caps the homeomorphism problem for closed 3-manifolds} is at least as hard as
{\caps the graph isomorphism problem}. The difficulty of the graph isomorphism problem
is far from clear. Among researchers, there are differing opinions about whether it might admit
a polynomial-time solution. However, note that Babai has recently given a remarkable algorithm that runs in quasipolynomial time [3].
A discussion of the computational complexity of the problem is given in Section 13.2 of [3], where Babai
says that `it is quite possible the intermediate status of GI [graph isomorphism] (neither NP-complete nor polynomial time)
will persist.' Theorem 1.2 therefore presents a potentially non-trivial lower bound
for the complexity of the homeomorphism problem for closed 3-manifolds.

Given Theorem 1.2, it is natural to ask whether the equivalence problem
for classical links is conditionally hard. We have not been able to establish such a result. Nevertheless, the
following related problem is shown to be conditionally hard.

\vskip 6pt
\noindent {\caps The sublink problem}
\vskip 6pt

A {\sl sublink} of a classical link $L$ is a union of components of $L$.
In this problem, one is given diagrams of two links $L_1$ and $L_2$, and one asks: is $L_1$ equivalent to a sublink of $L_2$?

We will prove the following result.

\noindent {\bf Theorem 1.3.} {\caps The sublink problem} {\sl is NP-hard.}

\vskip 6pt
\noindent {\caps Unresolved issues}
\vskip 6pt

It would, of course, be desirable to prove that many other problems on knots, links and 3-manifolds are conditionally hard.
For example, the problem of deciding whether two classical knots are equivalent 
{\sl appears} to be difficult. Certainly, all known solutions are highly intricate and appear not to admit efficient algorithmic
implementation. (See [19] for a survey.) However, one must draw a large distinction between problems that appear to be hard and
problems that are genuinely conditionally hard. All NP-complete problems are computationally equivalent to
{\caps sat}, the problem of deciding whether a collection of Boolean variables admits a truth-assignment
so that a given collection of sentences is satisfied. In other words, in order to show that a problem is conditionally hard, it needs to have
a sufficient level of generality or universality. Given this observation, it would appear to be difficult to show that
many problems in low-dimensional topology are conditionally hard. In this context, the results provided by Theorems 1.1, 1.2 and 
1.3 are perhaps rather surprising.

\vskip 6pt
\noindent {\caps Acknowledgement}
\vskip 6pt

The author would like to thank the referee for their suggestions which have undoubtedly improved this paper.

\vskip 18pt
\centerline {\caps 2. The Thurston complexity of an unoriented classical link}
\vskip 6pt

In this section, we prove Theorem 1.1. This is based on the proof by Agol, Hass and Thurston [2]
that the problem of deciding whether a knot in a compact orientable 3-manifold bounds a compact
orientable surface with genus at most $g$ is NP-complete. Their
argument has two main steps. In the first step, they show that the problem is in NP, and in the second
step, they show that it is NP-hard. Our proof of Theorem 1.1 also divides into these two steps.

We now explain why the problem {\caps Upper bound for the Thurston complexity of an unoriented classical link} is in NP. We are given a diagram
$D$ of a link $L$ and an integer $n$ (in binary). We suppose that $L$ bounds a Seifert surface $S$
with Thurston complexity at most $n$. We need to establish the existence of a certificate, which can be verified in polynomial
time as a function of the crossing number $c(D)$ and the number of digits of $n$, that proves that $L$ does indeed bound such a surface.

We use the diagram $D$ to construct a triangulation $T$ for the exterior $X = S^3 - {\rm int}(N(L))$. With a little care, one can arrange 
that its number of tetrahedra $t$ is at most a fixed linear function of the crossing number $c(D)$. 
This triangulation restricts to a triangulation of the boundary $\partial N(L)$, and we need to ensure that this
has a certain form, as described shortly. The surface $S$ induces an orientation on $L$. We include this orientation as part of our certificate.
Once this orientation has been specified, the curves $\partial S$ in $\partial X$ are determined up to isotopy. This is because
a curve on each torus is determined, up to isotopy, by the number of times it runs around the meridional and longitudinal
directions. In the case of a component of $\partial S$, the number of longitudes is $1$, and the number of meridians
is the linking number between the relevant component of $L$ and the remainder of the link. Thus, knowing only
the orientation on $L$, one knows the isotopy class of $\partial S$ in $\partial X$. Hence, we may arrange the triangulation
$T$ so that a representative for the isotopy class of $\partial S$ is disjoint from the vertices of $T$, transverse to the edges
and intersects each edge of $T$ at most once.
Even with this extra restriction, $T$ can still be constructed in polynomial time and its number of tetrahedra $t$
is still at most an explicit linear function of $c(D)$. Moreover the representative curves for $\partial S$ can
also be constructed explicitly in polynomial time.

The surface $S$ can be assumed to be incompressible. Hence, by a minor variation of Proposition 3.3.24 of [21], there is
a Seifert surface for $L$, that is homeomorphic to $S$ and with boundary equal to the specified curves
$\partial S$ in $\partial N(L)$, and that is in normal form. This means that it intersects each tetrahedron of $T$ in a
collection of normal triangles and squares (as in Figure 3.10 in [21]). We also call this normal surface $S$.

Associated with $S$, there is a list of $7t$ non-negative integers, which count the number of triangles
and squares of each type. This list is the {\sl normal surface vector} $[S]$. Normal surface vectors satisfy
a system of equations called the {\sl matching equations}. They also satisfy conditions called the
{\sl quadrilateral constraints} which assert that two different types of square cannot coexist
within the same tetrahedron (see Section 3.3.4 of [21] for example). Haken showed [11] that non-negative integer solutions
to the matching equations satisfying the quadrilateral constraints determine properly embedded normal
surfaces in $X$. 

A normal surface $S$ is {\sl fundamental} if $[S]$ cannot be
written as a sum $[S_1] + [S_2]$ for two other non-empty normal surfaces $S_1$ and $S_2$.
In the proof of Theorem 4.1.10 in [21], it is shown that, when $L$ has a single component
and the curve $\partial S$ intersects each edge of $T$ at most once (as we have already arranged),
then there is a Seifert surface $S$ of minimal Thurston complexity which is normal and fundamental and with boundary equal to this
specified curve. In the case where
$L$ has more than one component, this need not be the case, but we will show that $S$ can be
chosen to be a bounded sum of fundamental normal surfaces. More precisely, we will show
that $S$ can be chosen so that $[S] = [S_1] + \dots + [S_k]$, where each $S_i$ is fundamental
and $k$ is at most the number of components of $L$.

Choose $S$ so that it has minimal possible number of intersections with the 1-skeleton of $T$,
subject to the condition that it is normal,  has boundary the given curves in $\partial N(L)$
and satisfies $\chi_-(S) = \chi_-(L)$.

Now our condition that each edge of $T$ intersects $\partial S$ at most once implies that
if $S$ is a sum of normal surfaces $S_1$ and $S_2$, then each component of $\partial S_i$
is a component of $\partial S$.

Write $S$ as a sum of fundamental surfaces, and suppose that the number of these surfaces
exceeds the number of components of the link. Then, by the above observation, at least
one of these surfaces $S_i$ is closed. Note that it is connected because it is fundamental.
Let $S'_i$ denote the sum of the remaining surfaces.
Then $[S] = [S_i] + [S_i']$. Now it is shown in the proof of Theorem 4.1.10 of [21] that we may arrange that
no component of $S_i - {\rm int}(N(S'_i))$ or $S_i' - {\rm int}(N(S_i))$ is a disc in the interior of the manifold. 
(See also Lemma 2.1 in [16].)
Hence, no component of $S_i$ or $S_i'$ is a sphere. No component can be a
projective plane, because these do not embed in $S^3$. Also, if any component of
$S'_i$ is a disc, then this is disjoint from $S_i$, and so is a component of $S$.
Now, $S$ is obtained from the union of $S_i$ and $S_i'$ by smoothing off the intersection curves $S_i \cap S_i'$.
Hence, $\chi(S) = \chi(S_i) + \chi(S_i')$. 

Suppose that $S'_i$ is orientable. Then, writing $d(S)$ and $d(S'_i)$ for the number of
disc components of these surfaces,
$$\chi_-(S'_i) = -\chi(S'_i) + d(S'_i) = -\chi(S) + \chi(S_i) + d(S) \leq \chi_-(S).$$
The last inequality holds because $S_i$ is a closed connected surface other than a sphere or projective plane.
Thus, after discarding any closed components, $S'_i$ becomes a normal Seifert surface for $L$ with the same boundary as $S$, with no greater
Thurston complexity, but with fewer points of intersection with the 1-skeleton. This
is a contradiction.

Suppose now that $S'_i$ is not orientable. The components of $S_i - {\rm int}(N(S'_i))$ and $S_i' - {\rm int}(N(S_i))$
inherit an orientation from $S$. Suppose that we form the surface $S'$ by pasting together 
$S_i - {\rm int}(N(S'_i))$ and $S_i' - {\rm int}(N(S_i))$ in the manner that respects the
orientation. Then $S'$ is a Seifert surface for $L$ plus possibly some closed components. Its Euler characteristic is the
same as that of $S$, and it has the same number of disc components. It also has
no sphere components. Hence, $\chi_-(S') = \chi_-(S)$. But, as shown in Lemma 3.2.12 of [21],
there is an ambient isotopy of $S'$, leaving its boundary fixed, that reduces its number of
intersections with the 1-skeleton. Furthermore, it may then be transformed into a normal surface, without increasing
its Thurston complexity, without moving its boundary and without increasing its number
of intersections with the 1-skeleton. This is a contradiction. Thus, we have shown that $[S]$
is a sum $[S_1] + \dots + [S_k]$, where each $S_i$ is fundamental
and $k$ is at most the number of components of $L$.

Hass and Lagarias [13] gave an upper bound for the modulus of each co-ordinate for $[S_i]$ when $S_i$ is
fundamental. This was an exponential function of the number $t$ of tetrahedra in the given triangulation of $X$. 
Thus, when written in binary, the number of digits of this co-ordinate is bounded above by a linear function of
$t$. Hence, we obtain a similar upper bound for the number of digits of each co-ordinate of $[S]$.
This vector $[S]$ forms part of the certificate. More precisely, if there is a positive answer to the decision problem
{\caps Upper bound for the Thurston complexity of an unoriented classical link}, or, in other words, if there is a Seifert surface 
for $L$ with Thurston complexity at most $n$, then the certificate consists of the following data:
the orientation on $L$ induced by this Seifert surface and the vector $[S]$ for 
a normal Seifert surface $S$ which minimises Thurston complexity and with this exponential bound
on each co-ordinate of $[S]$. To verify the certificate, the algorithm of
Agol-Hass-Thurston [2] is applied. This produces a list of components for $S$, together with their Euler characteristic,
and so one can use this to verify that $\chi_-(S)$ is at most $n$. One also verifies that this surface has
boundary equal to the specified curves in $\partial N(L)$. This
completes our summary that the decision problem {\caps Upper bound for the Thurston complexity of an unoriented classical link}
is in NP.

The second part of the Agol-Hass-Thurston proof established that the problem of deciding whether a knot in a compact orientable 3-manifold
bounds a compact orientable surface with genus at most $g$ is NP-hard. To do this, they show that the decision problem {\caps 1-in-3-sat} is Karp-reducible to it. 
Here, one is given a collection of variables
$v_1, \dots, v_n$, and a collection of Boolean sentences. Each sentence is of the form `Exactly one of $x_1$, $x_2$ and $x_3$ is true',
where each $x_i$ is either a variable or its negation. We call this collection of variables and sentences an {\sl instance} of {\caps 1-in-3-sat}.
The problem asks whether there is an assignment of `true' or `false' to
each variable, so that all the sentences become true. We call this a {\sl solution} to this instance.

We will show that {\caps Upper bound for the Thurston complexity of an unoriented classical link} is NP-hard in a similar way, by establishing that a solution to
this decision problem leads to a solution of {\caps 1-in-3-sat}. More specifically, the latter problem is Karp-reducible to the
former one.
So, suppose that we are given an instance of {\caps 1-in-3-sat}, with variables $v_1, \dots, v_n$ and sentences $c_1, \dots, c_m$.
Each sentence is of the form `Exactly one of $x_i$, $x_j$ and $x_k$ is true'. We write this as $x_i \veebar x_j \veebar x_k$.
We now describe the construction of a diagram of an associated classical link $L$. We will show that $L$ has Thurston complexity at most $4m$ if
and only if the given instance of {\caps 1-in-3-sat} has a solution. It is the construction of $L$ where our
argument differs from that of Agol, Hass and Thurston [2] most substantially.

The diagram for $L$ will lie in the plane $\{ z = 0 \}$ in ${\Bbb R}^3$, and the diagrammatic projection map will be the usual
vertical projection map from ${\Bbb R}^3$ onto this plane. To form this diagram, 
we start with a bipartite graph embedded in ${\Bbb R}^3$. It has $n$ vertices spaced along the line $\{ y = 1, z = 0 \}$, corresponding
to the variables $v_1, \dots, v_n$, and $m$ vertices spaced along the line $\{ y = 0, z = 0 \}$, corresponding to the sentences $c_1, \dots, c_m$.
If a variable $v_i$ or its negation appears in the sentence $c_j$, then a nearly straight edge of the graph runs between the corresponding vertices.
These edges are only nearly straight, for two reasons. Firstly, we wish for the graph to be embedded in ${\Bbb R}^3$, and so crossings are inserted
to ensure that the edges do not intersect. Secondly, it is possible for the same variable to appear more than once within a sentence,
and so more than one edge may run between the same pair of vertices. We now add a vertex to this graph, called the
{\sl extra vertex}, on the line $\{ y = -1, z = 0 \}$,
and we join it to each of the vertices corresponding to the sentences, each via a single edge.

We now use this graph to form the diagram of the link $L$. We first of all choose a knot $K$ with $\chi_-(K)$ at least $2m+1$, 
and crossing number at most some polynomial function of $m$ and $n$ (where $m$ and $n$ are the number of sentences and variables, respectively).
For the sake of being precise, we set it to be the $(2, 2m+3)$ torus knot. We now replace each sentence vertex with a link consisting of
four components, as follows. We place a knotted solid torus within a regular neighbourhood of the vertex, and where the core curve of
this solid torus is a copy of $K$. We call this a {\sl sentence solid torus}.
We pick four simple closed curves on the boundary of this solid torus, each of which has winding number one around
the solid torus. So, this is a satellite link with companion $K$, and where the four components each run once around the satellite in the
simplest way possible. We abbreviate this diagrammatically as in Figure 1.

\vskip 12pt
\centerline{
\includegraphics[width=1in]{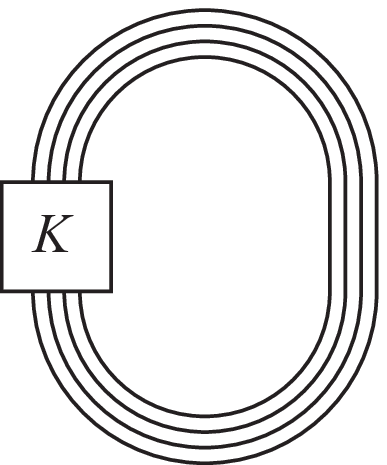}
}
\vskip 6pt
\centerline{Figure 1: The link of four components near each sentence vertex}

At each of the remaining vertices, we place a small round disc in the plane $\{ z = 0 \}$. Part of the boundary of this disc will form a subset of $L$.
When the disc is associated with a variable, we call it a {\sl variable disc}. The disc associated with the extra vertex is termed the
{\sl extra disc}.

Finally, we replace each edge of the graph by two parallel strands. Aside from the edges running to the extra vertex,
each edge corresponds to an occurrence of a variable $v_i$ or its negation $\neg v_i$ in a sentence $c_j$.
If the negation $\neg v_i$ occurs, then we add a half-twist to the two parallel strings. In the case of the edges running to the extra
vertex on the line $\{ y = -1, z = 0 \}$, we do not add half-twists.

We explain how to attach these parallel strands to the remainder of the link. One end runs to the boundary of the disc corresponding to $v_i$ or to the extra vertex.
The other endpoints of the two parallel
strings are banded onto a component of the four-component link corresponding to the sentence $c_j$. There are four edges coming 
into this sentence vertex, three corresponding to the variables in the sentence $c_j$, and the remaining one running from
the extra vertex. We ensure that the four bands are attached onto the four distinct components of the corresponding link.

\vfill\eject
Note that some arbitrary choices were made when constructing $L$, such as the crossing information between the edges
of the initial bipartite graph. These choices will not affect the main property of $L$ that we use, which is summarised in Proposition 2.1 below.
It is clear that the given diagram for $L$ has crossing number that is bounded above by some polynomial
function of $m$ and $n$.

An example of the link $L$ is shown in Figure 2. Here, there are three variables $v_1$, $v_2$ and $v_3$. There are three sentences
$$c_1 : v_1 \veebar \neg v_2 \veebar v_3, \qquad c_2 : v_1 \veebar \neg v_1 \veebar v_3, \qquad c_3 : \neg v_1 \veebar \neg v_2 \veebar \neg v_3.$$

Note that $L$ has $n+1$ components for the following reason. Each component is associated with a variable vertex or with the extra vertex.
It is obtained from the boundary circle of a variable disc or the extra disc by attaching `fingers'. Each finger runs along a band as far
as a sentence solid torus, once around this solid torus, and back along the band to the initial disc.

\vskip 12pt
\centerline{
\includegraphics[width=5in]{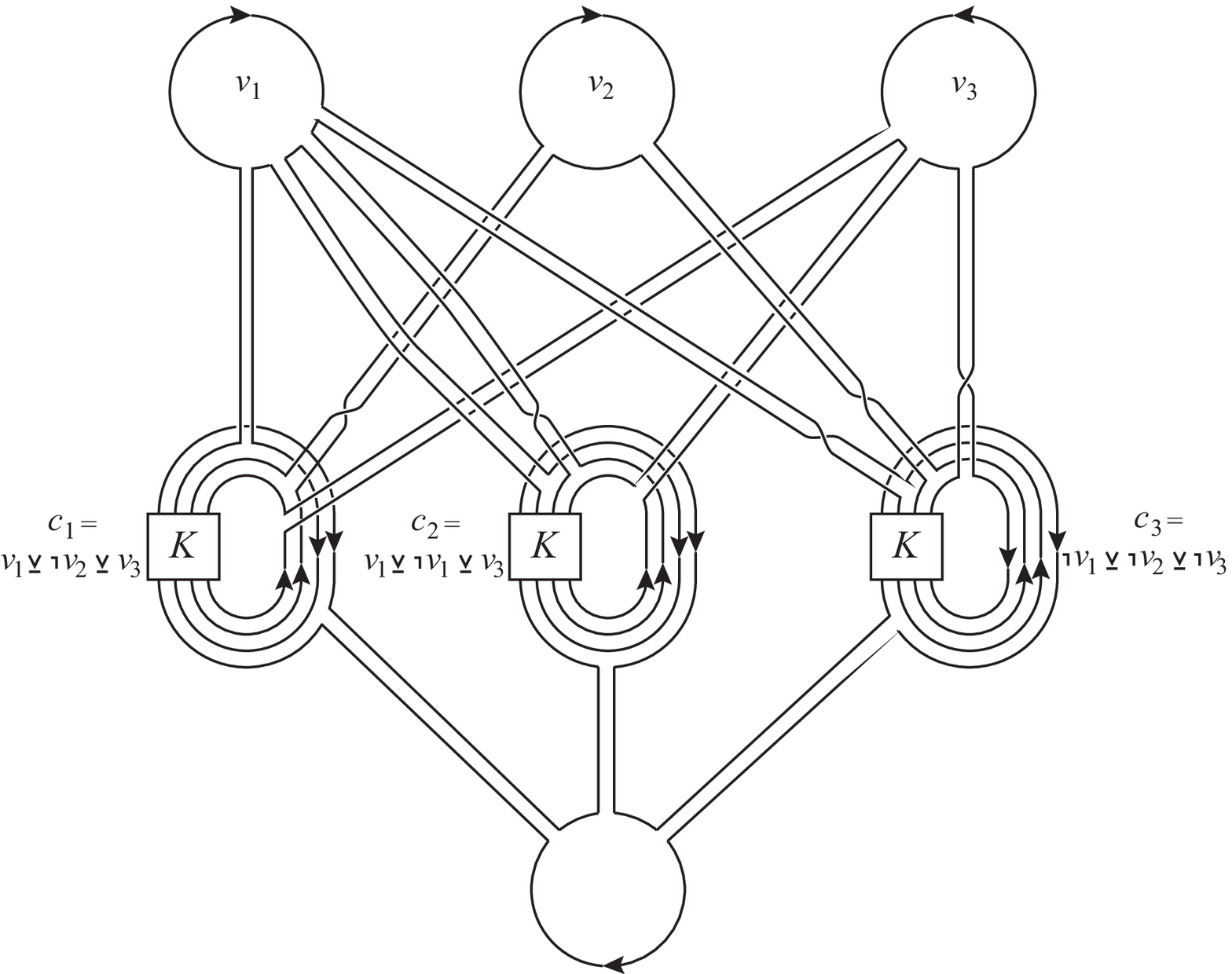}
}
\vskip 6pt
\centerline{Figure 2: The link $L$ associated with the variables $v_1$, $v_2$ and $v_3$ and sentences $c_1$, $c_2$ and $c_3$.} 
\centerline{A balanced orientation on $L$ is shown corresponding to the assignment $v_1 = T$, $v_2 = T$, $v_3 = F$.}

The following proposition will complete the proof of Theorem 1.1, because it shows that {\caps 1-in-3-sat} is Karp-reducible to 
{\caps Upper bound for the Thurston complexity of an unoriented classical link}. 

\noindent {\bf Proposition 2.1.} {\sl $L$ has Thurston complexity at most $4m$ if and only if the given instance of {\caps 1-in-3-sat} has a solution.}

One of the main tools in the proof of this is a well known lower bound on the Thurston complexity of satellite links. Recall that a link $L$
is a {\sl satellite} of a knot $K$ if it lies within a regular neighbourhood $N(K)$ of $K$. It is common to impose extra restrictions on the
way that $L$ lies within $N(K)$. For example, one often insists that $L$ does not lie within a 3-ball in $N(K)$ and that it is not
a core curve. We do not make these hypotheses here.

When $L$ is oriented, it then represents an element of $H_1(N(K))$.
This is $w$ times a generator, for some non-negative integer $w$, and we say that $w$ is the {\sl winding number} of this satellite.

For an oriented link $L$, we will only consider Seifert surfaces $S$ where the oriented boundary of $S$ matches the orientation
on $L$. Let $\chi_-^{\rm or}(L)$ denote the minimal Thurston complexity of a Seifert surface for $L$, subject to this condition.

The following lemma is well known. For example, a version appears as Proposition 2.10 in [5], where it is attributed to Schubert. 
We include a proof because we could not find a reference in print that deals with satellites in this level of generality.

\noindent {\bf Lemma 2.2.} {\sl Let an oriented link $L$ be a satellite of a knot $K$ with winding number $w$. Then
$\chi_-^{\rm or}(L) \geq w \ \chi_-(K)$.}

\noindent {\sl Proof.} Note first that we may assume that $K$ is non-trivial and that $w$ is non-zero, as otherwise the
statement of the lemma is empty.

Let $S$ be a Seifert surface for the oriented link $L$ with minimal Thurston complexity. This intersects the torus
$T = \partial N(K)$ in a collection of simple closed curves. If any of these curves is inessential in $T$, we may find one
that is innermost in $T$, and then compress $S$ along the disc that this bounds. The resulting surface may have
a closed component, but in this case, we discard this component. The resulting Seifert surface $S'$ satisfies $\chi_-(S') \leq \chi_-(S)$,
and so by the minimality of $\chi_-(S)$, we deduce that $\chi_-(S') = \chi_-(S)$. In this way, we may assume that
$S$ intersects $T$ in essential simple closed curves in $T$. 

We then deduce that each curve of $S \cap T$ is essential in $S$. For if there is an inessential curve, then an innermost
one in $S$ bounds a disc with interior disjoint from $T$. If this disc lies in $S^3 - {\rm int}(N(K))$, we deduce that
$K$ is the unknot. If the disc lies in $N(K)$, we deduce that $L$ lies within a 3-ball in $N(K)$ and hence that $w$ is zero.
We have assumed that neither of these possibilities arises. Since each curve of $S \cap T$ is essential in $S$,
we conclude that any disc components of $S$ are disjoint from $T$.

All the curves of $S \cap T$ are parallel in $T$. If
they are not all coherently oriented, then we may find two curves that are adjacent in $T$ and incoherently oriented.
We may cut $S$ along these two curves, and attach two copies of the annulus between them. The resulting
surface has the same Thurston complexity. Again, we discard any closed components. Since the new Seifert surface
intersects $T$ in fewer curves, we may assume that all curves of $S \cap T$ are coherently oriented.

The surface $S \cap N(K)$ forms a homology in $N(K)$ between $L$ and $S \cap T$. Hence, $S \cap T$ represents
$w$ times a generator of $H_1(N(K))$. The intersection between $S$ and $S^3 - {\rm int}(N(K))$ is trivialising
homology for these curves. Hence, $S \cap T$ forms $w$ longitudes on $\partial N(K)$. The intersection $S \cap (S^3 - {\rm int}(N(K)))$
forms a collection of $w$ Seifert surfaces for $K$, for the following reason. Suppose that some component $F$
of $S \cap (S^3 - {\rm int}(N(K)))$ had more than one boundary component. This component would either
be separating or non-separating in $S^3 - {\rm int}(N(K))$; we will show that both possibilities lead to a
contradiction. If $F$ is separating in $S^3 - {\rm int}(N(K))$, then every simple closed curve in $S^3 - {\rm int}(N(K))$
has zero algebraic intersection number with it. If $F$ is non-separating in $S^3 - {\rm int}(N(K))$, then
some simple closed curve has intersection number one with it. However, since $F$ has boundary equal
to $k > 1$ coherently oriented curves on $T$, its intersection number with a meridian $\mu$ for $K$
is a multiple of $k$. Since $[\mu]$ generates $H_1(S^3 - {\rm int}(N(K)))$, we deduce that
the intersection number between $F$ and every simple closed curve in $S^3 - {\rm int}(N(K))$ is
a multiple of $k$. This is incompatible with both the case where $F$ is separating and the
case where it is non-separating.

So, $S\cap T$ bounds $w$ Seifert surfaces for $K$.
Hence, $S \cap (S^3 - {\rm int}(N(K)))$ has Euler characteristic at most $- w \, \chi_-(K)$.
It is an essential subsurface of $S$. Therefore,
denoting the union of components of $S$ that are not discs by $S_-$, we have
$$\chi_-^{\rm or}(L) = -\chi(S_-) \geq -\chi(S \cap (S^3 - {\rm int}(N(K))) \geq w\,  \chi_-(K).$$ 
$\square$

\noindent {\sl Proof of Proposition 2.1.} We need to show that $L$ has Thurston complexity at most $4m$ if and only if the given instance of {\caps 1-in-3-sat} has a solution.
We will prove that these statements are both equivalent to the existence of a {\sl balanced} orientation on $L$,
which is defined as follows. An orientation on the components of $L$ is balanced if and only if, for
each sentence of the instance, two of the four components of the link within the associated solid torus are oriented one way and the other two are oriented the other.

Suppose first that there is a solution to this {\caps 1-in-3-sat} instance. Recall that $L$ has $n+1$ components. Each component
is obtained from the boundary circle of a variable disc or the extra disc by attaching fingers. So one may orient
this component simply by specifying an orientation of the boundary of this disc. For each variable $v_i$ given a true value, we orient the
boundary of the disc in a clockwise fashion around the disc. For each variable given a false value,
we orient the the boundary circle in an anti-clockwise fashion. For the disc corresponding to the extra vertex,
we orient its boundary clockwise. Consider the four link components that enter each sentence solid torus. 
One of these is attached to the
extra vertex; it is oriented clockwise. So is the string corresponding to the term $x_i$ which is true
within the sentence $x_i \veebar x_j \veebar x_k$. The other two strings are oriented anti-clockwise.
Hence, this is a balanced orientation on the link.

Conversely, suppose that $L$ has a balanced orientation. Then for each variable vertex and the extra vertex,
the strings that run around the associated disc are all oriented clockwise or all oriented anti-clockwise. This is because, as one encircles
the vertex, one runs around the vertex, then along a band to a sentence solid torus, then once around the sentence
solid torus, and then back along the same band. Hence, within a band, the two strings are incoherently oriented.
So, around a variable vertex, all the strings are coherently oriented.
 We may reverse the orientation of every component of $L$
and still get a balanced orientation. In this way, we can ensure that the strings of the extra vertex
are oriented clockwise. We give the variable $v_i$ the value true if and only if the strings that encircle the associated disc
are oriented clockwise. Now the strings encircling the extra vertex are oriented clockwise. Hence,
within each sentence solid torus, the string that is banded to the extra vertex is oriented clockwise.
Therefore, within that sentence solid torus, exactly one of the other strings is oriented clockwise,
and the remaining two are oriented anti-clockwise. So, within the associated sentence $x_i \veebar x_j \veebar x_k$,
exactly one of the statements $x_i$, $x_j$ and $x_k$ is true. In other words, the sentence $x_i \veebar x_j \veebar x_k$
is true. So, this is a solution to this {\caps 1-in-3-sat} instance.

Suppose now that we have a balanced orientation on $L$. Then we can construct an oriented Seifert surface for $L$
that induces this orientation on $L$, as follows. For each variable vertex and the extra vertex, insert the associated horizontal disc.
Orient this disc compatibly with the orientation on the strings in its boundary. For each sentence solid torus,
insert two oriented annuli into the solid torus, so that the boundary of these annuli equals the four strings
running on the boundary of the solid torus, and so that their orientations are correct. Now attach bands to these
discs and annuli, each running along an edge of the graph. The resulting surface is oriented, and its
oriented boundary gives the balanced orientation on $L$. It consists of $n+1$ discs and $2m$ annuli with $4m$ bands attached,
and so its Euler characteristic is
$n+1-4m$. The number of components of $L$ is $n+1$, and so this is an upper bound
for the number of disc components of the surface. So, this surface has Thurston complexity at most $4m$.

Conversely suppose that $L$ bounds a Seifert surface with Thurston complexity at most $4m$. 
This induces an orientation on $L$, which we claim is a balanced one. For suppose that in some sentence solid
torus, it is not the case that two of the strings are oriented one way and the other two are oriented the other.
Note that, by enlarging this solid torus, the entire link $L$ is a satellite of $K$. Because of our assumption
about the orientations of the strings, this satellite has winding number at least 2. So, by Lemma 2.2,
$\chi_-^{\rm or}(L) \geq 2 \chi_-(K) \geq 4m+2$, which is a contradiction. We deduce that this orientation on $L$ must
be balanced. $\square$

\vskip 18pt
\centerline{\caps 3. The Homeomorphism Problem for 3-manifolds}
\vskip 6pt

In this section, we prove that the {\caps graph isomorphism problem} is Karp-reducible to the {\caps homeomorphism problem for closed 3-manifolds}. 
The proof is little more than an observation, together with an application
of some basic 3-manifold theory. But it seems a worthwhile result nonetheless.

\noindent {\sl Proof of Theorem 1.2.} Suppose that we are given a solution to the homeomorphism problem. We use this to
provide a solution to the graph isomorphism problem. Let $\Gamma$ be a finite graph. From this, we construct
an associated triangulated orientable 
3-manifold $M(\Gamma)$. There are many possible ways of doing this, but here is one. For each vertex $v$ of $\Gamma$ with valence $d(v)$,
we form a triangulated copy of $P \times S^1$, where $P$ is a compact orientable surface with $d(v)$ boundary components and with genus $1$. 
We may choose this triangulation so that the number of tetrahedra is bounded above by a linear function of $d(v)$. We associate
each boundary component of $P \times S^1$ with an endpoint of an edge leaving $v$. In this torus, we specify two
slopes, the {\sl fibre} $\{ \ast \} \times S^1$ and the {\sl longitude} $C \times \{ \ast \}$, where $C$ is the relevant
boundary component of $P$. We orient these slopes using a fixed orientation of $S^1$ and of $P$.
We arrange that this torus is triangulated using just two triangles, and where the
fibre and longitude form edges. When two vertices are joined by an edge of the graph, we glue together
the corresponding tori via an orientation-reversing homeomorphism that swaps the longitude and the fibre and
that preserves their orientations. The result is the manifold $M(\Gamma)$.

We claim that two finite graphs $\Gamma_1$ and $\Gamma_2$ are isomorphic if and only if the 3-manifolds $M(\Gamma_1)$
and $M(\Gamma_2)$ are homeomorphic. For suppose that there is an isomorphism from $\Gamma_1$ to $\Gamma_2$.
We use this to build a homeomorphism from $M(\Gamma_1)$ to $M(\Gamma_2)$, as follows. If an
edge $e_1$ in $\Gamma_1$ is sent to an edge $e_2$ in $\Gamma_2$, then this homeomorphism
will send the torus in $M(\Gamma_1)$ corresponding to $e_1$ to the torus in $M(\Gamma_2)$ corresponding to $e_2$,
in such a way that longitudes and fibres are preserved. We need to show that this extends to a homeomorphism
from $M(\Gamma_1)$ to $M(\Gamma_2)$. Each complementary region of these tori in $M(\Gamma_1)$ corresponds to a vertex
of $\Gamma_1$. This is sent via the graph isomorphism to a vertex of $\Gamma_2$. These two vertices
have the same valence, and so they correspond to homeomorphic copies of $P \times S^1$. Moreover,
given any permutation of the boundary components of $P$, there is a homeomorphism of $P \times S^1$
that realises this. Hence, we may construct the required homeomorphism from the complementary
region in $M(\Gamma_1)$ to the complementary region in $M(\Gamma_2)$ that extends the
given homeomorphism on the tori.

Conversely, suppose that $M(\Gamma_1)$ and $M(\Gamma_2)$ are homeomorphic. Each manifold has
its canonical collection of JSJ tori. In the case of these graph manifolds, this is just the tori associated with
edges of the graphs, for the following reason. These tori divide the manifolds into Seifert fibre spaces.
None of these is a solid torus, because we assumed that the surfaces $P$ had genus $1$. Similarly, no two of
these tori are parallel, again because each surface $P$ had genus $1$. The Seifert fibrations on
adjacent pieces induce distinct fibrations of each torus. Hence, these tori form
the canonical tori, up to isotopy. (One can prove this by noting that they satisfy the definition of 
a `W-decomposition' as defined in [22], for example.) Therefore, the given homeomorphism $M(\Gamma_1) \rightarrow M(\Gamma_2)$
can be modified by an isotopy so that it sends one collection of tori to the other. This therefore
induces a bijection from the edges of $\Gamma_1$ to the edges of $\Gamma_2$. The homeomorphism 
sends each complementary region of one set of tori to a complementary region of the other. Hence, we obtain a bijection from
the vertices of $\Gamma_1$ to the vertices of $\Gamma_2$. These bijections between vertex sets and edge sets
respect the incidence of edges and vertices, and hence they specify an isomorphism $\Gamma_1 \rightarrow \Gamma_2$.
$\square$

\vskip 18pt
\centerline {\caps 4. The sublink problem}
\vskip 6pt

In this section, we will show that the {\caps sublink problem} is NP-hard by finding an NP-hard problem that Karp-reduces to it.
This is the problem of determining whether a finite graph contains a Hamiltonian path. Recall that a {\sl Hamiltonian path}
in a graph $\Gamma$ is a path that visits each vertex precisely once. Note that we are concerned here with paths
that start and end at distinct vertices, not paths which end where they start, thereby forming a cycle.

It is well known that the decision problem that asks whether a finite graph contains a Hamiltonian path is
NP-complete [9]. We will show that the {\caps sublink problem} is NP-hard, by showing that a solution to it
can be used to provide a solution to the {\caps Hamiltonian path problem}. Given a finite graph $\Gamma$, we wish to construct 
a diagram of a link $L(\Gamma)$. This will involve some arbitrary choices.
However, the construction will have the following property.

\noindent {\bf Proposition 4.1.} {\sl Let $\Gamma$ be a graph with $n$ vertices. Then $\Gamma$ contains
a Hamiltonian path if and only if $L(\Gamma)$ contains as a sublink a string of trefoils with $2n-1$ components.}

The {\sl string of trefoils} with $5$ components is shown in Figure 3. It is a collection of trefoils and unknots, arranged alternately
in a line. We require that at the start and end of this line, there are trefoils. This generalises in an obvious way to any odd
number of components.

\vskip 12pt
\centerline{
\includegraphics[width=3.5in]{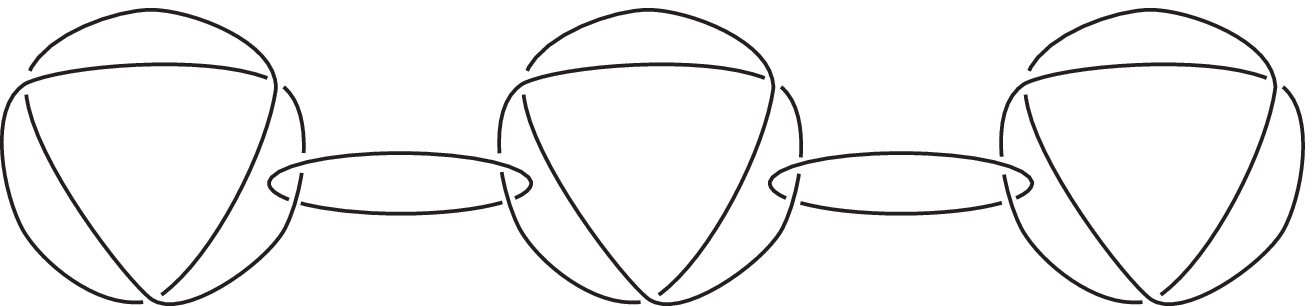}
}
\vskip 6pt
\centerline{Figure 3: A string of trefoils}

The diagram of $L(\Gamma)$ will contain at most $n^4+2$ crossings. Hence, given a solution to the
{\caps sublink problem}, we can apply it to $L(\Gamma)$, and this will determine whether $\Gamma$ contains
a Hamiltonian path.

The construction of $L(\Gamma)$ is as follows. If $\Gamma$ contains any edge loops, we remove them, since an edge loop
can never be part of a Hamiltonian path. Similarly, if any pair of vertices is joined by multiple edges, we replace these by 
a single edge. Let $\Gamma_-$ be the resulting graph. We next embed $\Gamma_-$ in ${\Bbb R}^3$, in a more-or-less
arbitrary way, as follows. Arrange the vertices of $\Gamma_-$ around a circle in the plane, and so that no three resulting diagonals intersect in a point. Whenever
two vertices are joined by an edge, realise this as a diagonal. These diagonals will intersect in general, and if so, we move one slightly off the
other, by introducing a crossing. The resulting diagram for $\Gamma_-$ has at most $n^4/8$ crossings, since any pair of edges cross at most once,
and there are at most $n(n-1)/2 < n^2/2$ edges.

We now use this diagram of $\Gamma_-$ to create a diagram for $L(\Gamma)$. Replace each vertex by a diagram of a small trefoil in the plane
with three crossings. Replace each edge by an unknotted curve, which runs along the edge, around the trefoil at one endpoint, back along the edge, and then
around the trefoil at its other endpoint. We may choose this curve so that it projects to an embedded circle in the plane.
Whenever two edges of $\Gamma_-$ have a crossing, this gives rise to four crossings of $L(\Gamma)$. Also,
each edge of $\Gamma_-$ gives rise to two crossings at each of its endpoints. Hence, the diagram of $L(\Gamma)$ has
at most $n^4/2 + 2n^2 + 3n$ crossings. This is less than $n^4+2$ as long as $n \geq 3$, and it is easy to check that
$n^4+2$ is also a bound for the crossing number when $n =1$ or $2$. In Figure 4, we give an example of
$L(\Gamma)$ in the case where $\Gamma$ is the complete graph on $5$ vertices.

We now prove Proposition 4.1. It is clear that if $\Gamma$ contains a Hamiltonian path, then the corresponding vertices and
edges of this path gives rise to a sublink of $L(\Gamma)$ that is a string of trefoils. This has $2n-1$ components, since the path contains
$n$ vertices and $n-1$ edges. An example is shown in Figure 4, where
a Hamiltonian path in $\Gamma$ induces the sublink shown in bold. 

Conversely, suppose that $L(\Gamma)$ contains a sublink that is a string of trefoils with $2n-1$ components. We note that the components of
$L(\Gamma)$ have the following properties. The linking number of any two components of $L(\Gamma)$ has modulus either zero or one.
Moreover, it is one if and only if one of the components corresponds to an edge, the other corresponds to a vertex and
this vertex is incident to this edge. Hence, as one travels along the sublink that is a string of trefoils, the corresponding components
of the sublink correspond to an alternating sequence of vertices and edges of the graph. Moreover, successive components
of the sublink correspond to a vertex and an edge that are incident. So, this corresponds to a path in $\Gamma$,
where all the vertices and edges of the path are distinct. There are $n$ distinct vertices in this path, and hence
this forms a Hamiltonian path in $\Gamma$. $\square$

\vskip 12pt
\centerline{
\includegraphics[width=3in]{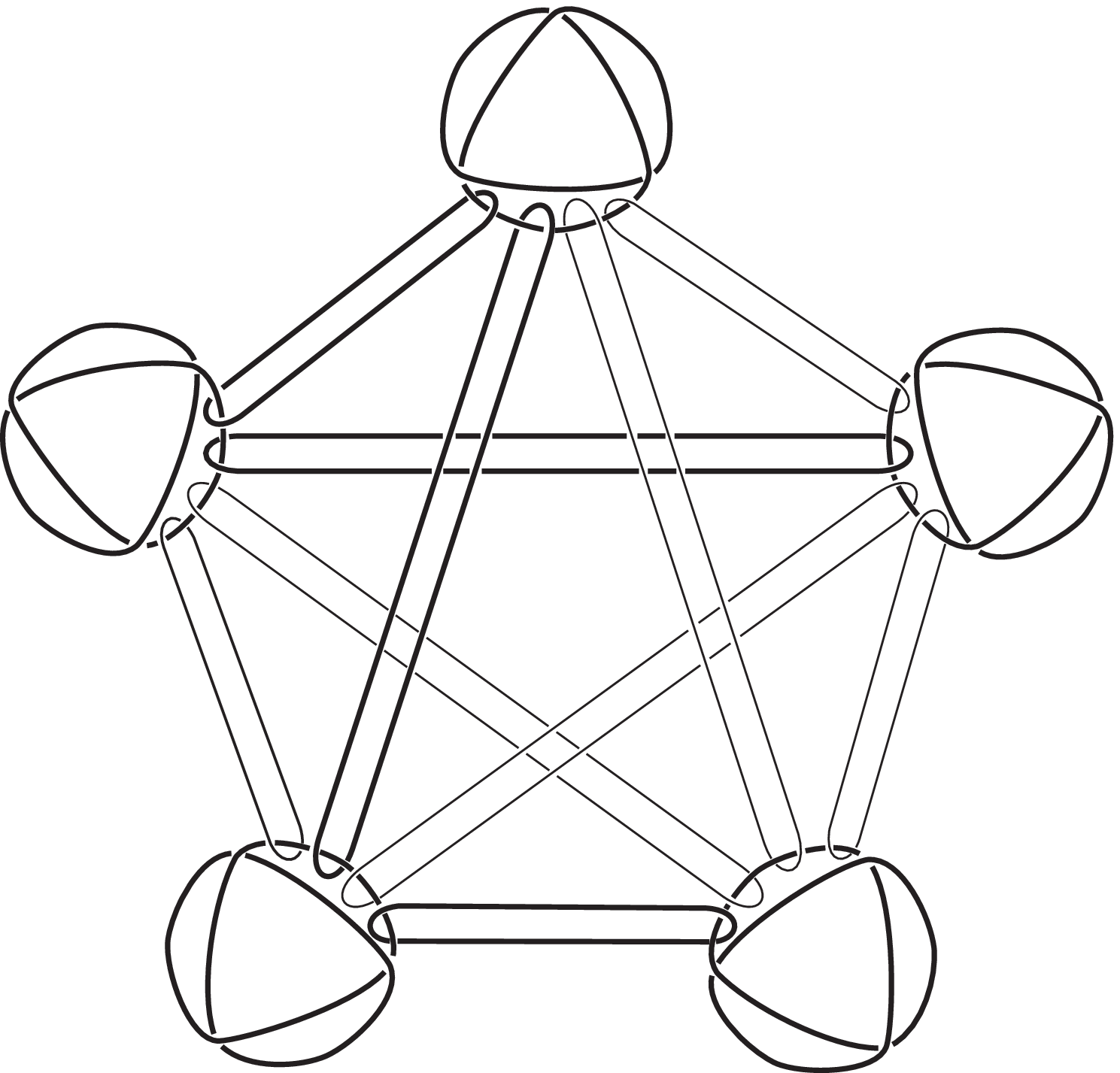}
}
\vskip 6pt
\centerline{Figure 4: The link $L(\Gamma)$ arising from the complete graph $\Gamma$ on five vertices.}
\centerline{The sublink corresponding to a Hamiltonian path in $\Gamma$ is shown in bold.}

Note that we deliberately chose to focus on the {\caps Hamiltonian path problem}, rather than the more
familiar NP-complete problem of deciding whether a finite graph contains a Hamiltonian cycle. This
is because a Hamiltonian cycle in $\Gamma$ could possibly form a knotted simple closed curve
in the given embedding, and so the corresponding sublink of $L(\Gamma)$ would not be of
a specific link type. Moreover, unknotted Hamiltonian cycles need not specify a unique sublink type in $L(\Gamma)$,
because they might have varying writhes. By using Hamiltonian paths rather than cycles, we avoid these issues.

We have shown that the {\caps sublink problem} is NP-hard. We do not know whether it is likely to be NP-complete.
However, it is certainly algorithmically solvable using the known solution to the equivalence problem for links.

\vskip 18pt
\centerline{\caps 5. Conclusion}
\vskip 6pt

It seems likely that many problems in low-dimensional topology are not solvable in polynomial time. Certainly,
this is the case for the three problems considered in this paper, provided $P \not= NP$ and provided the {\caps graph
isomorphism problem} has no polynomial-time solution. It remains a significant challenge to show that other problems
are conditionally hard. Possibly the most notable outstanding case is the equivalence problem for knots and links.
The goal is to reduce it to a problem lying outside of low-dimensional topology and that is widely
viewed as unsolvable in polynomial time.

\vskip 18pt
\centerline {\caps References}
\vskip 6pt

\item{1.} {\caps I. Agol,} {\sl Knot genus is NP}, Conference presentation (2002).
\item{2.} {\caps I. Agol, J. Hass, W. Thurston}, {\sl The computational complexity of knot genus and spanning area},
Trans. Amer. Math. Soc. 358 (2006) 3821--3850.
\item{3.} {\caps L. Babai,} {\sl Graph Isomorphism in Quasipolynomial Time,} arXiv:1512.03547
\item{4.} {\caps L. Bessi\`eres, G. Besson, S. Maillot, M. Boileau, J. Porti,} {\sl Geometrisation of 3--manifolds,}
EMS Tracts in Math. 13, European Math. Soc. (EMS), ZŸrich (2010)
\item{5.} {\caps G. Burde, H. Zieschang,} {\sl Knots,} de Gruyter Studies in Mathematics, 5. Walter de Gruyter \& Co., Berlin, 2003.
\item{6.} {\caps B. Burton, E. Colin de Verdi\`ere, A. de Mesmay}, {\sl On the complexity of immersed normal surfaces},
arXiv:1412.4988.
\item{7.} {\caps B. Burton, A. de Mesmay, U. Wagner}, {\sl Finding non-orientable surfaces in 3-manifolds}, 
arXiv:1602.07907
\item{8.} {\caps B. Burton, J. Spreer}, {\sl The complexity of detecting taut angle structures on triangulations.}
Proceedings of the Twenty-Fourth Annual ACM-SIAM Symposium on Discrete Algorithms, 168Ð183, SIAM, Philadelphia, PA, 2012. 
\item{9.} {\caps M. Garey, D. S. Johnson} {\sl Computers and Intractability: A Guide to the Theory of NP-Completeness,} W.H. Freeman (1979).
\item{10.} {\caps O. Goldreich}, {\sl Computational complexity. A conceptual perspective}, Cambridge University Press (2008).
\item{11.} {\caps W. Haken,} {\sl Theorie der Normalfl\"achen.} Acta Math. 105 (1961) 245--375.
\vfill\eject
\item{12.} {\caps W. Haken,} {\sl Some results on surfaces in $3$-manifolds.}
Studies in Modern Topology pp. 39--98 Math. Assoc. Amer. (1968)
\item{13.} {\caps J. Hass, J. Lagarias,} {\sl The number of Reidemeister moves needed for unknotting.}
J. Amer. Math. Soc. 14 (2001), no. 2, 399--428 
\item{14.} {\caps J. Hass, J. Lagarias, N. Pippenger,} {\sl The computational complexity of knot and 
link problems.} J. ACM 46 (1999), no. 2, 185--211.
\item{15.} {\caps G. Hemion}, {\sl  On the classification of homeomorphisms of $2$-manifolds 
and the classification of $3$-manifolds.} Acta Math. 142 (1979), no. 1-2, 123--155. 
\item{16.} {\caps W. Jaco, U. Oertel,} {\sl An algorithm to decide if a $3$-manifold is a Haken manifold}, 
Topology 23 (1984) 195--209.
\item{17.} {\caps G. Kuperberg}, {\sl Knottedness is in NP, modulo GRH}, Adv. Math. 256 (2014), 493--506.
\item{18.} {\caps G. Kuperberg}, {\sl Algorithmic homeomorphism of 3-manifolds as a corollary of geometrization}, 
arXiv:1508.06720.
\item{19.} {\caps M. Lackenby}, {\sl Elementary knot theory}, To be published by the Clay Mathematics Institute.
\item{20.} {\caps M. Lackenby}, {\sl The efficient certification of knottedness and Thurston norm}, arXiv:1604.00290
\item{21.} {\caps S. Matveev}, {\sl Algorithmic topology and classification of 3-manifold},
Algorithms and Computation in Mathematics, 9. Springer, Berlin, 2007.
\item{22.} {\caps W. Neumann, G. Swarup}, {\sl Canonical decompositions of 3-manifolds},
Geom. Topol. 1 (1997) 21--40.
\item{23.} {\caps G. Perelman,} {\sl The entropy formula for the 
Ricci flow and its geometric applications,} Preprint,
arxiv:math.DG/0211159
\item{24.} {\caps G. Perelman,} {\sl  Ricci flow with surgery on three-manifolds,}
Preprint, arxiv:math.DG/0303109
\item{25.} {\caps G. Perelman,} {\sl Finite extinction time for the 
solutions to the Ricci flow on certain three-manifolds,} Preprint,
arxiv:math.DG/0307245 
\item{26.} {\caps P. Scott, H. Short}, {\sl The homeomorphism problem for closed 3-manifolds.}  
Algebr. Geom. Topol. 14 (2014), no. 4, 2431--2444.

\vskip 12pt
\+ Mathematical Institute, University of Oxford, \cr
\+ Radcliffe Observatory Quarter, Woodstock Road, Oxford OX2 6GG, United Kingdom. \cr

\end